\documentclass{article}

\usepackage[a4paper,outer=3cm]{geometry}

\usepackage[utf8]{inputenc}
\usepackage[T1]{fontenc}
\usepackage{times}

\usepackage{xspace}

\usepackage{longtable}
\usepackage{amsmath,amsthm,amssymb,amsthm}
\usepackage{enumerate}

\usepackage{algorithm}
\usepackage{algpseudocode}
\algnewcommand\algorithmicinput{\textbf{Input:}}
\algnewcommand\algorithmicoutput{\textbf{Output:}}
\algnewcommand\Input{\item[\algorithmicinput]}
\algnewcommand\Output{\item[\algorithmicoutput]}

\usepackage{hyperref}

\usepackage{caption}
\usepackage{subfigure}
\usepackage{tikz}
\usepackage{svg}

\usepackage{pdftricks}
\begin{psinputs}
	\usepackage{pstricks,pst-plot,pst-node,pst-func}
\end{psinputs}
\usepackage{graphics,graphicx}

\tikzstyle{vertex}=[circle, draw, inner sep=0pt, minimum size=6pt]

 \usepackage{marginnote}

\newcommand{\m}[1]{}

\usepackage[nothm]{thmbox}
\usepackage{amsthm}
\usepackage{thmtools}

\declaretheorem[parent=section,thmbox=M]{theorem}

\declaretheorem[numberlike=theorem,thmbox=M]{proposition}
\declaretheorem[numberlike=theorem,thmbox=M]{corollary}
\declaretheorem[numberlike=theorem,thmbox=M]{conjecture}

\declaretheorem[numberlike=theorem]{lemma}
\declaretheorem[numberlike=theorem]{remark}
\declaretheorem[numberlike=theorem]{claim}

\newcommand{\Gya}{Gy\'arf\'as\xspace}

\newcommand{\Chu}{Chudnovsky\xspace}

\newcommand{\Hag}{H\"aggkvist\xspace}

\newcommand{\mF}{\mathcal{F}}

\newcommand{\Ra}{\Rightarrow}

\newcommand{\ora}[1]{\overrightarrow{#1}}

\newcommand{\ovlra}{\overleftrightarrow}

\DeclareMathOperator{\dic}{\ora \chi}

\newenvironment{maintheorem}{%
    \medskip
  \thmbox[M]{\textbf{Theorem \ref{thm:struct}}}%
  \hspace*{-1.5em}\slshape\ignorespaces%
  }
  {%
  \endthmbox\vspace*{.75ex}%
  }
  
\newenvironment{inroundtheorem}{%
    \medskip
  \thmbox[M]{\textbf{Theorem \ref{thm:in-round}}}%
  \hspace*{-1.5em}\slshape\ignorespaces%
  }
  {%
  \endthmbox\vspace*{.75ex}%
}

\newcommand{\SD}{S^+_2}

\title{Decomposing and colouring some locally semicomplete digraphs}
\author{Pierre Aboulker$^1$, Guillaume Aubian$^{1,2}$, Pierre Charbit$^{2}$\\
\small ($1$) DIENS, \'Ecole normale sup\'erieure, CNRS, PSL University, Paris, France\\
\small ($2$) Université de Paris, CNRS, IRIF, F-75006, Paris, France.}

\begin{document}

\maketitle
\begin{abstract}
    A digraph is semicomplete if any two vertices are connected by at least one arc and is locally semicomplete if the out-neighbourhood and the in-neighbourhood of any vertex induce a semicomplete digraph. In this paper we study various subclasses of locally semicomplete digraphs for which we give structural decomposition theorems. As a consequence we obtain several applications, among which an answer to a conjecture of Naserasr and the first and third authors: if an oriented graph is such that the out-neighbourhood of every vertex induces a transitive tournament, then one can partition its vertex set into two acyclic digraphs.
\end{abstract}

\section{Introduction}


In this paper all directed graphs (\textit{digraphs} in short) are simple, i.e. contain no loop and no multi-arc. If in addition a digraph contains no \textit{digon} (a cycle on two vertices), we say that it is an {\em oriented graph}.

For a vertex $x$ of a digraph $D$, we denote by $x^+(D)$ (resp. $x^-(D)$) the set of its out-neighbours (resp. in-neighbours) :
\begin{eqnarray*}
    x^+(D)=\{y\in V(D),\ xy\in A(D)\}\\
    x^-(D)=\{y\in V(D),\ yx\in A(D)\}
\end{eqnarray*}
If there is no ambiguity on the digraph, we will simply use $x^+$ and $x^-$.

Given a digraph $D=(V,A)$, its {\em underlying graph} is the unoriented graph $G=(V,E)$ such that $xy\in E$ if $xy\in A$ or $yx\in A$. We say that a digraph is {\em connected} if its underlying graph is connected as an unoriented graph. A {\em semicomplete} digraph is a digraph whose underlying graph is a complete graph. A {\em tournament} is a semicomplete oriented graph. 
A digraph $D$ is \textit{strongly connected} if there is a directed path between each ordered pair of vertices. In this case we say that $D$ is a \textit{strong} digraph. 

We denote by $\ovlra{K_2}$ the digraph on two vertices with an arc in both directions, by  $C_k$ the directed cycle on $k$ vertices and by $TT_k$  the unique acyclic tournament on $n$ vertices, called {\em transitive tournament}.
The oriented graph on three vertices with a vertex of out-degree $2$ (resp of in-degree $2$) and two vertices of in-degree $1$ (resp. two vertices of out-degree $1$) is called $S^+_2$ (resp. $S^-_2$).

For a class $P$ of digraphs (like semicomplete, tournament, acyclic), a digraph is {\em locally out-P} (resp. {\em locally in-P}) if for every vertex $x$, $x^+$ (resp. $x^{-}$) induces a digraph in $P$. 
For example, a digraph $D$ is locally out-semicomplete if the out-neighborhood of every vertex of $D$ induces a semicomplete digraph.
Finally, we will say that a digraph is {\em locally P} if it is both locally out-P and locally in-P. 
We make one exception for one of the main class of this paper: for the oriented graphs for which the out-neighborhood of every vertex is a transitive tournament, we will use the term "out-transitive oriented graphs" instead of the heavier and possibly confusing "out-transitive tournament oriented graphs".


A {\em linear order} on a digraph $D$ is an order $O:=(v_1,v_2,\ldots,v_n)$ of its vertices. Two orders $O_1$ and $O_2$ are equivalent if $O_1:=(v_1,v_2,\ldots,v_n)$ and $O_2=(v_k,v_{k+1},\ldots,v_n,v_1,v_2,\ldots, v_{k-1})$ for some $k$. An equivalence class for this relation is called a {\em cyclic order} of $D$. 
For two vertices $v_i$ and $v_j$ in a linear order the \textit{cyclic interval} $[v_i,v_j]$ is defined as follow:
$$
 [v_i,v_j] =     
\begin{cases}
     \{v_k, k\in [i,j]\} \text{ if $i<j$}\\
     \{v_k, k\not\in ]j,i[ \} \text{ if $i\geq j$}
    \end{cases}
$$ 
Note that cyclic intervals only depend on the cyclic order and not on a linear order chosen as a representative. As usual, $]v_i,v_j[ = [v_i, v_j] \setminus \{v_i, v_j\}$.
\medskip

Semicomplete digraphs are well studied and a natural and fruitful way to extend results on this class is to look at the class of locally semicomplete digraphs.
Introduced in 1990 by Bang-Jensen \cite{B90}, locally semicomplete digraphs have since then been the topic of more than 100 research papers, and a whole chapter in \cite{BGCD} is devoted to this class. A particularly nice result in this area is one of Huang that gives a geometric characterization of locally transitive digraphs. We state it here in the particular case of oriented graphs.

\begin{theorem}[Huang, \cite{H92}]\label{thm:round}
If $D=(V,A)$ is a connected oriented graph, then the two conditions below are equivalent 
\begin{enumerate}
    \item for every vertex $x$, both $x^+$ and $x^-$ induce a transitive tournament. 
    \item there exists a cyclic order of the vertices of $D$ such that 
$$\forall xy\in A, \forall z\in ]x,y[, xz\in A \text{ and } zy\in A$$
\end{enumerate}
\end{theorem}

Any oriented graph (strong or not) that satisfies the second condition above  is called a \textit{round} oriented graph. In other words, for every vertex $x$, $x^+$ (resp. $x^-$) consists in a cyclic interval placed just after (resp. before) $x$ in the cyclic order. Note that a round oriented graph is strong if and only if every vertex has at least one in-neighbour. 
This is because the cyclic order given by the theorem then yields an Hamiltonian cycle. By the similar observation, if a round oriented graph is not strong, then it is in fact acyclic.
\medskip

Our first result is a generalization of the theorem above in the particular case of strong oriented graphs. 

\begin{figure}[htbp]
\begin{center}
\includegraphics[width=0.25\textwidth]{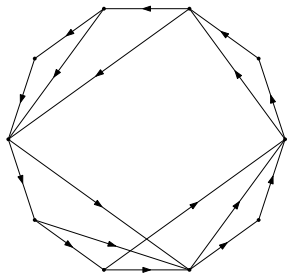}
\caption{An in-round oriented graph that is not round}
\label{fig:inround}
\end{center}
\end{figure}

\begin{inroundtheorem}
Let $D$ be a strong oriented graph. Then conditions below are equivalent.
\begin{enumerate}
    
    \item \label{itm:in-round_fst} for every vertex $x$, $x^+$ induces a tournament and $x^-$ induces an acyclic digraph 
    \item \label{itm:in-round_snd} there exists a cyclic order of the vertices of $D$ such that 
$$\forall xy\in A, \forall z\in ]x,y[, zy\in A$$

\end{enumerate}
\end{inroundtheorem}
Again condition \ref{itm:in-round_snd} can be seen as the property that for every vertex $x$, $x^-$ consists in a cyclic interval placed just before $x$ in the order (see Figure \ref{fig:inround}). Following the terminology of \cite{BGCD}, any oriented graphs satisfying condition \ref{itm:in-round_snd} of the theorem above will be called {\em in-round}. 

In \cite{LZM08}, we note that the authors prove a similar result giving an alternate condition 1: they ask that for every vertex $x$, $x^+$ induces a {\em transitive} tournament and $x^-$ induces an acyclic digraph {\em with a hamiltonian path}. These additional conditions (which are easily seen in fact to be implied by condition 2) are unnecessary, which makes our theorem slightly stronger. Moreover, our proof, exposed in section~\ref{sec:dec}, is also much shorter.

In fact, Theorem \ref{thm:in-round} is the first step towards our main structural result, which is a decomposition theorem for the class of strong locally out-transitive digraphs. If $H$ is a subdigraph of $D$, we define the {\em contraction} $D/H$ as the digraph obtained by removing all vertices of $H$, then adding a new vertex $h$ such that $xh$ (resp. $hx$) is an arc of $D/H$ if $x^{+}\cap H$ (resp $x^{-}\cap H$) is non empty. Beware that in general $D/H$ might contain digons (even if $D$ does not), eventhough in our case this never happens.

\begin{maintheorem}
If $D$ is a strong locally out-transitive oriented graph, then there exists a partition of its set of vertices into strong subdigraphs $D_{1},\ldots,D_{k}$ such that the digraph obtained by contraction of every $D_{i}$ is a strong in-round oriented graph.
\end{maintheorem}

We give two applications of this theorem. The first one is an acyclic colouring result: we prove that if $D$ is an out-transitive oriented graph, then it can be partitioned into two acyclic induce subdigraphs. This proves a conjecture of~\cite{ACN21}. Proof together with context of the conjecture is given in Section~\ref{sec:colouring}. This result was independently proof by Raphael Steiner, see Remark~\ref{Steiner}. 
The second one is a proof that locally in-tournaments satisfy the famous Caccetta-\Hag conjecture. This result is mentioned in  \cite{R13} as an unpublished but not easy result of Paul Seymour but here the idea is to show how our decomposition  theorem more or less directly implies the result. The proof is given in Section~\ref{sec:ch}

Finally, in Section \ref{sec:semicomplete}, we use the techniques developed for the proof of Theorem \ref{thm:main} to prove a structural theorem (see Theorem \ref{thm:k12k21}) for the class of locally-semicomplete digraphs. Such a theorem is given in \cite{BGGV} (see Theorem 3.12), but we feel that our theorem is simpler and might have some interesting applications. 
We propose some illustrations after the proof of the theorem.

One idea behind the paper is to promote the idea of finding  decomposition theorem for classes of digraphs, that is a theorem whose statement is of the kind: either a graph in this class is "basic" (belongs to some simple subclass), or it can be decomposed in some prescribed ways. These kinds of decomposition theorems proved to be a very powerful strategy in the world of undirected graphs (the most famous example being the celebrated proof of the perfect graph conjecture by \Chu et al in \cite{Perfect}), and there are to our knowledge not so many theorems of this kind in the world of directed graphs, and there is no reason to believe that it could not be as effective in this setting. 

\begin{remark}\label{Steiner}
A week prior to the submission of this paper, R. Steiner published on arXiv a very nice paper \cite{RS21} containing another proof of the colouring result of locally out-transitive oriented graphs mentioned above (as well as other results about acyclic colourings). Even though some of the ingredients are in common, the proof is longer and is different as it is an entirely inductive proof whereas ours relies on the structure theorem mentioned above.
\end{remark}


\section{Decomposing locally out-transitive oriented graphs}\label{sec:dec}

We start by proving the theorem mentioned in the introduction about strong oriented graphs that are both locally out-tournament and locally in-acyclic. 
\begin{theorem}\label{thm:in-round}
Let $D$ be a strong oriented graph. Then conditions below are equivalent.
\begin{enumerate}
    
    \item  for every vertex $x$, $x^+$ induces a tournament and $x^-$ induces an acyclic digraph 
    \item  there exists a cyclic order of the vertices of $D$ such that 
$$\forall xy\in A, \forall z\in ]x,y[, zy\in A$$

\end{enumerate}
\end{theorem}
We recall that we will use the term in-round for any oriented graph satisfying condition 2.

\begin{proof}[of Theorem \ref{thm:in-round}]
The easy direction is \ref{itm:in-round_snd} implies \ref{itm:in-round_fst}. Indeed, let $x,y,z$ be such that ${y,z}\subset x^+$ and assume w.l.o.g that $z\in [x,y]$. Then by \ref{itm:in-round_snd} we have that $zy\in A$, so $x^+$ is a tournament. 
Assume now that $x^-$  contains a directed cycle $C$, and let $y$ be the vertex of $C$ such that $C\setminus{y} \subset [y,x]$ (the left most vertex of $C$ in the representative of the cyclic order which ends in $x$). Let $z$ be the predecessor of $y$ on $C$. Now we have that $x\in [zy]$ and so by \ref{itm:in-round_snd} there must be an arc $xy$, which contradicts the fact that $y\in x^-$ (there are no digon here).

Now assume \ref{itm:in-round_fst}. For every vertex $x$, $x^-$ induces a non empty acyclic oriented graph, and hence contains a vertex $y$ such that $y^+\cap x^- = \emptyset$ (take the last vertex in a topological ordering of $x^-$). 
For every $x$ we arbitrarily choose one such vertex and denote it by $f(x)$. If $z$ is an out-neighbour of $f(x)$, then since $f(x)^+$ induces a tournament, $z$ and $x$ must be connected by an arc, and this cannot be $zx$ by definition of $f(x)$, so there must be an arc $xz$. We have therefore $f(x)^+ \setminus \{x\} \subset x^+$ for all $x$.

Now let $H$ be the graph induced by the arcs $f(x)x$. Each vertex of $H$ has in-degree exactly $1$, so $H$ contains a cycle $C$. If $C$ does not span all vertices of $D$, then since $D$ is strong, there exists an arc $xy$ in $D$ where $x\in C$ and $y\not\in C$. Because of the property above the whole cycle $C$ must be contained in $y^-$, which contradicts \ref{itm:in-round_fst}. So $H$ consists in an Hamiltonian cycle and the property $f(x)^+\subset x^+$ exactly translates into \ref{itm:in-round_snd} for the cyclic order defined by this Hamiltonian cycle.
 
\end{proof}

\begin{remark}\label{rem:transitive}
We point out that the proof above shows that Property \ref{itm:in-round_snd} implies in fact that every out-neighbourhood in a in-round oriented graph induces a {\em transitive} tournament (even if it is not necessary to add it in Property \ref{itm:in-round_fst} to get the equivalence). Thus in-round strong oriented graphs can also be seen as strong oriented graphs that are locally out-transitive and locally in-acyclic as proved proved in~\cite{LZM08}. 
\end{remark}


%

We now recall the statement of the main structural theorem mentioned in the introduction.

\begin{theorem}\label{thm:struct}
If $D$ is a strong locally out-transitive oriented graph, then there exists a partition of its set of vertices into strong subdigraphs $D_{1},\ldots,D_{k}$ such that the digraph obtained by contraction of every $D_{i}$ is an strong in-round oriented graph.
\end{theorem}


Let us start with a lemma about locally out-semicomplete digraphs (be aware that the rest of this section is about oriented graphs and not digraphs, but we state it here for digraphs as we will make use of it in Section~\ref{sec:semicomplete}, where we study digraphs).

\begin{lemma}\label{lem:base}
Let $D$ be a locally out-semicomplete digraph. Let $H$ be a strong subdigraph of $D$ and let $z\notin H$. If $z^-\cap H \neq \emptyset$ and $z^+\cap H = \emptyset$, then $H \subseteq z^-$.
\end{lemma}
\begin{proof}
Let $h$ and $h'$ be vertices of $H$ such that $hz$ and $hh'$ are arcs of $D$. Because $D$ is out-semicomplete and $z^+\cap H = \emptyset$, $h'z$ must be an arc. Since $H$ is strong, we get that $h^-\cap H = H$.
\end{proof}

\begin{proof}[of Theorem ~\ref{thm:struct} ]

From now on $D$ will denote a strong  locally out-transitive oriented graph. 

We define a {\em hub}  to be a subset of vertices $H$ of $D$ such that 
\begin{itemize}
\item $H$ induces a strong oriented graph,
\item there exists $x\notin H$ such that $H\subseteq x^{-}$.
\end{itemize}

Note that a hub is necessarily a strict subset of $V(D)$. 
A hub is trivial if it is reduced to one vertex. 

Recall that an oriented graph is in-round if the in-neighborhood of each vertex induces an acyclic oriented graph, and the out-neighborhood induces a tournament. Hence, since a cycle that is out-dominated is a hub, $D$ is in-round if and only if there is no non trivial hub. 

Assume $D$ is not in-round, and consider a hub $H$ that is (inclusion-wise) maximal. We want to prove that $D/H$ is a locally out-transitive oriented graph. Assume first by contradiction that there exists a vertex $z\notin H$ that is {\em mixed} for $H$, that is such that there exists $\{h,h'\}\subset H$ with $(hz,zh')\in A(D)^{2}$. Since $H$ is strong, there is a directed path from $h$ to $h'$. Due to the fact that $h^+$ induces a  tournament, by considering such $h$ and $h$' with a path of minimal length, one can assume that $hh'$ is an arc. There exists $x$ in $D$ such that $H\subset x^{-}$, and because $D$ does not contain $\SD$ there must be an arc between $x$ and $z$. It cannot be $xz$ because in that case $h^+$ would contain the directed triangle $xzh'$ which is not possible so $zx$ must be an arc. But then $H\cup\{z\}$ is a hub that contradicts the maximality of $H$.

So using Lemma \ref{lem:base}, we observe that eventually there are three types of vertices outside $H$: the one that have no arc to or from $H$, the one that have out-neighbours in $H$ but no in-neighbour, and the ones that are out-neighbour of every vertex of $H$. One can now easily check that $D/H$ 
is locally out-transitive.
\medskip

Consider now two distinct maximal hubs $H_{1}$ and $H_{2}$. We want to prove that their intersection is empty. For $i=1,2$, let $h_1$ be a vertex that in-dominates $H_i$. Note that it is not possible that $h_1\in H_2$ and $h_2\in H_1$ simultaneously, as this would imply a digon between the two vertices. So without loss of generality, assume $h_1\notin H_2$. But now if $H_1\cap H_2 $ is non empty, then $h_1$ has a in neighbour in $H_2$ and by Lemma \ref{lem:base} and since we proved that no vertex can be mixed for $H_2$, we have that $h_1$ in-dominates and $H_{1}\cup H_{2}$, which is therefore also a hub (it is strongly connected because the intersection is non empty). This contradicts the maximality of $H_1$ and $H_2$ 

Eventually we have proven that maximal hubs define a partition of the vertex set (it is trivial that every singleton is a hub). Moreover, by what precedes, if there is an arc $xy$ from a maximal hub $H_{1}$ to a maximal hub $H_{2}$, then $y$ is an out-neighbour of every vertex in $H_{1}$. We summarize this with the following claim (the transitive tournament fact is due to the graph being locally out-transitive).

\begin{claim} Maximal hubs form a partition of the vertex set. Moreover, if $H_{1}$ and $H_{2}$ are two maximal hubs, then between them:
\begin{itemize}
\item either there is no arc,
\item either all arcs from $H_{1}$ to a subset of $H_{2}$ inducing a transitive tournament and no other arc,
\item either all arcs from $H_{2}$ to a subset of $H_{1}$ inducing a transitive tournament and no other arc.
\end{itemize}

\end{claim}

Let $D'$ be the digraph obtained by contracting every maximal hub. If $D'$ contains a non trivial hub $H'$, then by the fact above it is clear that the set $H$ of vertices of $D$ that are mapped to vertices in $H'$ by the contraction form a hub that would contradict the maximality of the hubs that were contracted to obtain $D'$. So $D'$ contains no non trivial hub, so no dominated cycle, and hence by Theorem \ref{thm:in-round}, it is in-round. This concludes the proof of Theorem \ref{thm:struct}.
\end{proof}

\section{Applications of Theorem \ref{thm:struct}}
\subsection{Acyclic Colourings}\label{sec:colouring}
An {\em acyclic colouring} of a digraph is an assignment of colours to the vertices such that each colour class induces an {\em acyclic subdigraph}, that is a subdigraph containing no directed cycle.  The \textit{acyclic chromatic number}, or simply \emph{dichromatic number}, of a digraph $D$, is defined to be the smallest number of colours required for an acyclic colouring of $D$ and is dentoed $\dic(D)$. This notion was first introduced in 1982 by Neumann-Lara \cite{NL82} and has attracted  a lot of attention in the past decade (see for example \cite{ACL19, BHL18,  HK15, H17,  KS20, LM17, GSS20, S20}) 
as it seems to be the natural generalization for digraphs of the usual chromatic number. 

A \textit{hero} is a tournament $H$ such that every tournament not containing $H$ have bounded dichromatic number. In \cite{hero}, Berger et al give a full description of heros. 
Extending the notion of hero, a set $\mathcal{F}$ of digraphs is said to be  \emph{heroic} if every digraph 
not containing any member of $\mF$ as an induced subdigraph 
has bounded dichromatic number. 
In \cite{ACN21}, the first and third authors together with Naserasr proposed the following conjecture, along with a proof of the only if part. It can be seen as an analogue of the well-known Gy\'arf\'as-Sumner Conjecture  \cite{Gya87, Sum81} for oriented graphs.  

\begin{conjecture}[\cite{ACN21}]\label{conj:oriented}
Let $H$ be a hero and let $F$ be an oriented forest. The set $\{\ovlra{K_2}, H, F\}$ is heroic if and only if:
\begin{itemize}
\item either $F$ is the disjoint union of oriented stars,
\item or $H$ is a transitive tournament.
\end{itemize}
\end{conjecture}

This conjecture is still widely open. The first case that was left in \cite{ACN21} is the case $F=S^+_2$, and $H=C_{3}$ and it was conjectured that digraphs 
without any induced subgraph in $\{\ovlra{K_2}, C_3, S^+_2\}$
 have dichromatic number at most two. 
As mentioned in the introduction, we now prove this result, and in fact a stronger result. We recall that in~\cite{RS21}, Raphael Steiner independently proved the same result as well as other special cases of Conjecture~\ref{conj:oriented}, see Remark~\ref{Steiner}.  



\begin{theorem}\label{thm:2dic}
Every locally out-transitive oriented graph has dichromatic number at most $2$.
\end{theorem}

Note that it indeed extends the question mentioned above as it amounts to forbidding $\ovlra{K_2}$, $S^+_2$ and the tournament on $4$ vertices built by taking a directed triangle $C_{3}$ and adding a vertex with an arc going to the three other vertices. As already said, forbidding $S^+_2$ implies that every out-neighbourhood induces a tournament, and forbidding this 4-vertex tournament implies that every out neighbourhood must induce a $C_3$-free tournament, hence acyclic (every tournament containing a directed cycle must contain a $C_3$).

To ease the induction proof we will prove the following stronger result. 

\begin{theorem}\label{thm:main}
Let $D$ be a locally out-transitive oriented graph and $T$ a subdigraph of $D$ inducing a transitive tournament. Then $D$ admits a proper $2$-dicolouring such that all vertices of $T$ receive the same colour.
\end{theorem}

The first step towards the proof of this theorem is to prove the exact same statement for the class of in-round oriented graphs (which we recall can also be seen as locally out-transitive and locally in-acyclic oriented graphs).
\begin{proposition}\label{prop:colin-round}
Every in-round oriented graph has dichromatic number at most $2$. More precisely, for every vertex $x$, there exists a valid $2$-dicolouring such that $\{x\}\cup x^{+}$ is monochromatic.
\end{proposition}
\begin{proof}
We only need to prove it when the oriented graph is strong since a $2$-colouring of each strong component yields a valid $2$-colouring of the whole oriented graph. Consider the cyclic order given by the definition of in-round and pick any vertex $x$. Let $y$ be the vertex such that $xy$ is a longest arc, that is the arc such that the interval $[x,y]$ contains the maximum number of vertices. This implies that in the linear order given by the interval $]y,x[$  all arcs go forward since a back arc $x'y'$ would force the arc $xy'$ contradicting the maximality of $xy$. So $]y,x[$ induces an acyclic oriented graph. Moreover $[x,y]$ induces an acyclic oriented graph since it is included in $y^{-}$ and by definition of the in-round cyclic order. This concludes the proof.
\end{proof}

We are now ready to prove Theorem \ref{thm:main}.

\begin{proof}[of Theorem \ref{thm:main}]
Again we can assume that our digraph $D$ is strong, and we proceed by induction on the number of vertices of $D$. Let $T$ be a transitive tournament in $D$. We consider the decomposition into maximal hubs given by Theorem \ref{thm:struct} and label the hubs $H_{1},\ldots,H_{k}$ so that the order corresponds to the cyclic order $h_{1},\ldots,h_{k}$ of the contracted in-round oriented $D'$. Since the theorem guarantees that $D'$ contains no digons, there can exist arcs in only one direction between two distinct $H_i$. Moreover, Lemma \ref{lem:base} implies that if there is an arc $xy$ from $H_{i}$ to $H_{j}$, then there is an arc from every vertex in $H_{i}$ to $y$ and by the in-round property of $D'$ that the same is true for any vertex in $H_k$ for $k\in [i,j[$. We define $T_{i}$ for $i=1, \ldots, k$ to be the set of vertices in $H_{i}$ that have an in-neighbour out of $H_{i}$. By the remark just above, it in fact consists in $x^{+}\cap H_{i}$ for every $x\in H_{i-1}$. Note also that if a cycle of $D$ is not included in a maximal hub, then if it intersects $H_{i}$, it must intersect $T_{i}$. Let $s$ be the source in the transitive tournament $T$, and without loss of generality assume it belongs to $H_{1}$. Note that because every other vertex in $T$ is an out-neighbour of $s$, we have that $T\cap H_{i}\subset T_{i}$ for every $i\geq 2$.

Now we are ready to define our colouring. First we consider the in-round oriented graph $D'$ and denote by $T'$ the subtournament (transitive) in $D'$ consisting in all $h_{i}$ such that $T\cap H_{i}$ is non empty. By Proposition \ref{prop:colin-round}, there exists a valid $2$-dicolouring of $D'$ such that the out-neighbourhood of $h_{1}$ (which contains $T'$) is monochromatic, say coloured $1$. Now by induction we can ask for every $i\geq 2$ for a colouring of $H_{i}$ such that all vertices of $T_{i}$ get the colour of $h_i$ in the colouring of $D'$ defined above. For $i=1$ we ask by induction for a colouring of $H_{1}$ such that every vertex of $T\cap H_{1}$ gets the colour of $h_1$, that is $1$.

First, note that in this colouring, every vertex of $T$ gets colour $1$ because of the assumption on the colouring of $D'$. We now need to prove that that this colouring is a valid $2$-dicolouring of $D$. Assume by contradiction that $C$ is a monochromatic cycle and consider a minimal such cycle. $C$ is not included in some $H_{i}$ since the colouring is valid in every hub. So $C$ hits more than one hub, and since the vertices in each $H_{i}$ have the same out-neighbours out of $H_{i}$, the minimality of $C$ implies that $C$ contains at most one vertex from each hub $H_{i}$ and this vertex must belong to $T_{i}$ (otherwise it cannot have an in-neighbour out of $H_{i}$). The only case where this does not yield directly a monochromatic cycle (and hence a contradiction) in the contracted digraph $D'$ is if $C$ is coloured $2$ and contains a vertex $x$ in $T_{1}\setminus T$. Let $y$ be its successor on the cycle. Then $y$ belongs to some $T_{j}$ where $j$ is such that $h_{1}h_{j}$ is an arc of $D'$. But by assumption on the colouring of $D'$ this implies that $h_{j}$ gets colour $1$ and therefore $y$ must get colour $1$, which is a contradiction that finishes the proof of Theorem \ref{thm:main}.

\end{proof}


\subsection{About Caccetta-\Hag Conjecture}\label{sec:ch}
A beautiful and famous conjecture due to Caccetta and \Hag states states the following. 
\begin{conjecture}[Caccetta-\Hag]
Let $k\geq 2$ be an integer. Every digraph $D$ on $n$ vertices with no directed circuits of length at most $k$ contains a vertex of out-degree less than $n/k$.
\end{conjecture}
The case $k=2$ is trivial but the case $k=3$ is still widely open and has attracted a lot of attention. In \cite{R13} (see page 3), it is mentioned that for $k=3$, while adding the hypothesis that the graph has no $S^+_2$ makes it very easy, the dual case of forbidding $S^-_2$ was proven by Seymour but is "substantially more difficult". 

Here we prove that this comes as an easy consequence of Theorem~\ref{thm:struct} and Theorem~\ref{thm:in-round}, for any value $k\geq 3$.

\begin{theorem}\label{thm:CH}
Let $D$ be locally in-tournament oriented graph on $n$ vertices with no directed cycle of length at most $k$. Then $D$ contains a vertex of out-degree less than $n/k$.
\end{theorem}

Theorem~\ref{thm:struct} and Theorem~\ref{thm:in-round} were designed for locally out-tournament but of course by reversing the arcs we get an equivalent statement for locally in-tournament (here an out-round digraph is a digraph obtained by an in-round digraph by reversing all the arcs). We combine both to get this corollary.

\begin{corollary}\label{cor:locain}
If $D$ is a  
strong locally in-tournament that does not contain any $C_3$, 
there exists a partition of its set of vertices into strong subdigraphs $D_{1},\ldots,D_{k}$ such that the digraph $D'$ obtained by contraction of every $D_{i}$ is a strong oriented graph which is out-round, that is admits a cyclic order on its vertices such that 
$$\forall x,y,z\in V(D')\  (xy\in A \land  z\in ]x,y[  )\ \Ra\ xz\in A
$$
\end{corollary}

\begin{proof}[of Theorem \ref{thm:CH}]
Let $D$ be digraph as in the statement of the theorem and denote by $n$ the number of its vertices. First observe that we can assume $D$ to be strong, since one can apply it to a terminal strong component. Since $k\geq 3$, notice that $D$ does not contain $C_3$ so we can apply Corollary \ref{cor:locain} above. 

We begin by proving a wighted version for out-round oriented graphs.

\begin{lemma}
Let $D$ be a strong out-round oriented graph with no cycle $C_3$ and let $w$ be a positive weight function on the vertices of $D$. Denote by $W$ the sum of all weights. Then there exists a vertex $u$ such that 
$$\sum_{v\in u^+} w(v) < \frac{W-w(v)}{k} $$
\end{lemma}

Let us show first that this lemma indeed implies the theorem. Since $D$ admits an out-round quotient by Corollary \ref{cor:locain}, we can apply the lemma to the quotient with weights being the sizes of the contracted subdigraphs. Let $D_i$ be the subdigraph corresponding to the vertex of small out-degree given by the claim. Now by applying induction to $D_i$, one finds a vertex in $D_i$ with out-degree (strictly) less than $|D_i|/k$, which combined with the lemma gives the desired result.

Let us now prove the lemma. Consider the cycle order as in Corollary~\ref{cor:locain}. For every vertex $x$, denote by $f(x)$ the last of its out-neighbour along the cyclic order and by $\phi(x)$ the quantity $\sum_{y\in x^+} w(y)$. Observe that due to the out-round structure, $\phi(x)$ is exactly the sum of weights in the interval $]x,f(x)]$. We denote also by $f^{(i)}(x)=f(f^{(i-1)}(x)$ with the convention $f^{(0)}(x)=x$.

Consider a path on $k$ vertices of the form $xf(x)f(f(x))\ldots f^{(k-1)}(x)$. Observe that because of the out-round structure and the non existence of a cycle of length at most $k$ the path does not cycle completely around the cyclic order. Furthermore for the same reason we have that any two intervals $]f^{(i)}(x)f^{(i+1)}(x)]$ for distinct values $i\neq i'$ are disjoint. Eventually we have that for every vertex $x$
$$\sum_{i=0}^{k-1} \phi^{(i)}(x) \leq W-w(x)$$
Now consider some cycle $C$ formed by the arcs $xf(x)$. By summing the above inequality for all vertices in the cycle one gets
$$k\sum_{x\in V(C)} \phi(x) \leq W|C|-\sum_{x\in V(C)} w(x)$$
So there must exist a vertex in $V(C)$ that satisfies :
$$ k\phi(x)+w(x)\leq W $$
which is exactly the assertion of the lemma.

\end{proof}



%
%
%


\section{Structure of locally semicomplete digraphs} \label{sec:semicomplete}

In this section digons are allowed and we focus on locally semicomplete digraphs which we recall is the class of digraphs such that the in-neighbourhood and out-neighbourhood of any vertex is semicomplete.

A statement in the flavour of Conjecture \ref{conj:oriented} is easy to prove. The proof was independently obtained by Raphael Steiner~\cite{RS21}. 

\begin{theorem}
For every hero $H$, $H$-free locally tournaments have bounded dichromatic number. More precisely, it is bounded by at most twice the maximum dichromatic number of a $H$-free tournament. 
\end{theorem}

\begin{proof}[]
Let $k$ be the maximum dichromatic number of an $H$-free tournament. 
Let $x$ be any vertex. By induction on the number of vertices, we know that the oriented graph induced by the set of non neighbours $N=V\setminus (\{x\}\cup x^+ \cup x^-)$ is $2k$ colourable so we first properly colour $N$ with colours in $\{1,\ldots,2k\}$. Now since $x^+$ and $x^-$ are tournaments, they are $k$-colourable by hypothesis. We can thus use colours $\{1,\ldots,k\}$ for a proper colouring of $x^-$, and colours $\{k+1,\ldots,2k\}$ for a proper colouring of $x^+$, give $x$ any colour, and it is not difficult to check that this gives a valid $2k$-colouring of $D$ (the main fact to observe is that there are no arcs from $x^-$ to $N$, and from $N$ to $x^+$).
\end{proof}

The purpose of this section is to use our techniques to give a proof of the following structural theorem for the class of locally semicomplete digraphs. As already mentioned in the introduction, such a theorem is given in \cite{BGGV} (see Theorem 3.12), but our theorem is much simpler to state and might have some interesting applications. In fact the theorem of \cite{BGGV} has the same first two cases below, but their third case (called {\em evil} in Chapter 6 of the monograph \cite{BG01}) has a more complicated description that seem to us less easy to handle for applications. We propose some illustrations after the proof of the theorem.

We need some vocabulary first. A vertex $x$ {\em out-dominates} (resp. {\em in-dominates}) a set of vertices $X$ if $X\subseteq x^+$ (resp. $X\subseteq x^-$). A vertex $x$ {\em strictly out-dominates} (resp. {\em strictly in-dominates}) a set of vertices $X$ if $X \subseteq x^+\setminus x^-$ (resp. $X \subseteq x^- \setminus x^+$). A set of vertices $X$ {\em out-dominates} (resp. {\em strictly out-dominates}, {\em in-dominates}, {\em strictly in-dominates}) a set of vertices $Y$ if every vertex of $X$ out-dominates (resp. {\em strictly out-dominates}, {\em in-dominates}, {\em strictly in-dominates}) $Y$. Also, we call {\em universal vertex} a vertex $x$ such that $x^+=x^-=V(D)\setminus \{x\}$. 

\begin{theorem}\label{thm:k12k21}
Let $D$ a connected locally semicomplete digraph, then either :
\begin{itemize}
    \item $D$ is semicomplete with a universal vertex.
    \item There exists a partition of $V(D)$ into $k\geq 2$ subsets each inducing strong connected semicomplete digraphs such that the digraph obtained by contracting every member of the partition is a round oriented graph.
    \item there exists a partition of $V(D)$ into four sets $E$, $F$, $G$ and $H$ such that :
        \begin{itemize}
            \item $F$ and $H$ are non empty, and one of $E$ and $G$ is non empty.
            \item $D[E]$, $D[F]$, $D[G]$ and $D[H]$ are semicomplete;
            \item $E$ strictly out-dominates $F$, $F$ strictly out-dominates $G$, $G$ out-dominates $H$ and $H$ out-dominates $E$.
            \item $\forall x\in G,\ x^+\cap E\neq\emptyset$  and $x^-\cap E\neq\emptyset$
        \end{itemize}
\end{itemize}
\end{theorem}

\begin{figure}[htbp]
\begin{center}
\includegraphics[width=0.3\textwidth]{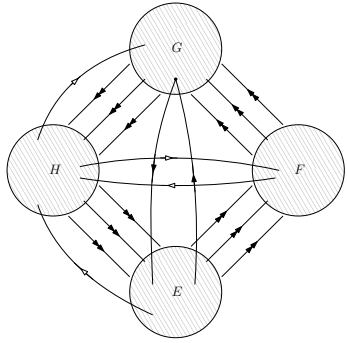}
\caption{The structure in the third case of Theorem \ref{thm:k12k21}}
\label{fig:4blobs}
\end{center}
\end{figure}

For this proof, we need to relax the notion of hub of the previous section: a set $X \subseteq V$ is called a {\em weak hub} if there exists a vertex which strictly in-dominates {\em or} strictly out-dominates $X$. The "strictly" part of this definition was implicit before for hubs since we were in the context of oriented graphs. A weak hub $X$ is said to be {\em mixed} if there exist vertices $x \notin X$ and $u,v \in X$ such that both  $xu$, $vx$ are arcs.  

We split the proof of Theorem \ref{thm:k12k21} into two parts depending on the existence of a maximal weak hub that is mixed. When there are no such subsets, the proof will have a lot of similarities with the one of Theorem \ref{thm:struct}.
\begin{lemma}\label{lem:no_mixed_weak hub}
Let $D$ be a connected locally semicomplete digraph such that no maximal weak hub is mixed. Then either $D$ is semicomplete with a universal vertex, or there exist $k \geq 2$ disjoint sets $X_1, X_2 \dots X_k$ of vertices such that $D[X_i]$ is strong for $i=1, \dots, k$, and  the digraph obtained by contraction of every $X_i$ is a round oriented graph. 
\end{lemma}

\begin{proof}[of Lemma \ref{lem:no_mixed_weak hub}]
Let $D$ be a connected locally semicomplete digraph and assume it is not semicomplete with a universal vertex. We first prove that maximal weak hubs form a partition of $V(D)$.

\begin{claim}\label{clm:weak hub_intersection}
Let $X$ and $Y$ two distinct maximal weak hubs. Then $X \cap Y = \emptyset$.
\end{claim}

{\bf Proof :}
Assume $X$ and $Y$ contradicts the claim. If there exists $x \notin Y$ which strictly in-dominates $X$, then there is at least one arc to $x$ from a vertex of $X \cap Y \neq \emptyset$ and since $Y$ is not mixed, there is no arc from $x$ to $Y$. As $D[Y]$ is strongly connected, Lemma \ref{lem:base} implies that $x$ strictly in-dominates $Y$. But then $X \cup Y$ would be a weak hub, contradicting the maximality of $X$.
If now there exists $x\in Y$ which strictly in-dominates $X$, then consider any arc $ab$ where $a\in X \cap Y$ and $b\in X\setminus Y$, which must exist since $D[X]$ is strong. But now $Y$ is mixed because of the arcs $ab$ and $bx$, a contradiction. 
So $X$ is not strictly in-dominated and is similarly it cannot be strictly out-dominated. 
\qed
\medskip
 
\begin{claim}\label{clm:weakhubcover}
Every vertex belongs to some maximal weak hub. 
\end{claim}
{\bf Proof :}
Let $u$ be a vertex of $D$. It suffices to prove  $\{u\}$ is a weak hub. Assume by contradiction this is not the case. By definition of weak hub this forces $u^+=u^-$. If $u^{+} = V \setminus \{u\}$, $D$ is semicomplete with a universal vertex, a contradiction. If not, as $D$ is connected, there exists vertices $v \in u^{+}$ and $w \notin u^{+}$ such that $vw \in A$ or $wv \in A$. But then, as $u^{-} = u^{+}$, $u$ and $w$ are both  in-neighbours or both out-neighbours of $v$,  implying $w \in u^{-} = u^{+}$, a contradiction. 
\qed

The next claim describes the structure of the arcs linking maximal weak hubs.
\begin{claim}\label{clm:weakhubinteract}
Let $X$ and $Y$ be two distinct maximal weak hubs. Then either $X$ strictly in-dominates $Y$, or $X$ strictly out-dominates $Y$, or there is no arc between $X$ and $Y$.
\end{claim}
{\bf Proof :}
Assume that there exists $x$ in $X$ and $y$ in $Y$ such that $xy\in A(D)$. By Lemma \ref{lem:base} applied to $H=X$ and $z=y$, and because $X$ is not mixed by hypothesis, we have that $y$ strictly in-dominates $X$. But now by applying the same Lemma to the digraph obtained from $D$ by reversing all arcs (which is still locally out-semicomplete since $D$ is locally semicomplete), we get that $Y$ strictly in-dominates $X$.
\qed

Let us now consider a partition of $V$ into maximal weak hubs $X_{1}, X_{2} \dots X_{k}$ and let $D'$ the digraph obtained by contracting every set in the partition. By Claim \ref{clm:weakhubinteract}, $D'$ is a locally tournament oriented graph. Moreover, the in-neighbourhood (resp. the out-neighbourhood) of any vertex $x \in V(D')$ is acylic. Indeed, if there were any such cycle $C$ on vertices of $D'$  corresponding to weak hubs $X_{i_{1}}, X_{i_{2}} \dots X_{i_{\ell}}$ with $\ell \geq 2 $, then $\cup_{j=1}^{\ell} X_{i_j}$ would induce a weak hub in-dominated (resp. out-dominated) by every vertex of the weak hub corresponding to $x$ in $D$, a contradiction to the maximality of the $X_i$. Due to Theorem \ref{thm:round}, $D'$ is a round digraph. 
\end{proof}

To prove Theorem \ref{thm:k12k21} it remains to prove the following lemma which deals with the case where there exists a maximal weak hub that is mixed.
\begin{lemma}\label{lem:mixed_weak hub}
Let $D$ a connected locally semicomplete digraph such that there exists a maximal weak hub which is mixed. Then there exists a partition of $V(D)$ into four sets $E$, $F$, $G$ and $H$ such that :
        \begin{itemize}
            \item $F$ and $H$ are non empty, and one of $E$ and $G$ is non empty.
            \item $D[E]$, $D[F]$, $D[G]$ and $D[H]$ are semicomplete;
            \item $E$ strictly out-dominates $F$, $F$ strictly out-dominates $G$, $G$ out-dominates $H$ and $H$ out-dominates $E$.
                        \item $\forall x\in G,\ x^+\cap E\neq\emptyset$  and $x^-\cap E\neq\emptyset$

        \end{itemize}
\end{lemma}

\begin{proof}[of Lemma \ref{lem:mixed_weak hub}]
Let $X$ be a maximal weak hub which is mixed. Let us define :
\begin{eqnarray*}
X^{M} &=& \{u \notin X \ \mid \exists a, b \in X \text{ such that } ua \in A, bu \in A \} \\
X^{IN} &=& \{u \notin X \ \mid \exists a \in X \text{ such that } ua \in A \}\setminus X^M\\
X^{OUT} &=& \{u \notin X\ \mid \exists a \in X \text{ such that } au \in A \}\setminus X^M \\
X^{NO} &=&V \setminus (X \cup X^{IN} \cup X^{OUT} \cup X^{M}) = \{u \in V \mid \forall v \in X, uv \notin A\text{ and } vu \notin A\}. 
\end{eqnarray*}

Since $X$ is strong, by Lemma \ref{lem:base}, $X$ strictly out-dominates $X^{OUT}$  and due to the same lemma applied to the digraph obtained by reversing all arcs of $D$, we have that $X^{IN}$ strictly out-dominates $X$. 

Eventually since $X$ is assumed to be a weak hub that is mixed, we have that $X^{M}$ and $X$ are non-empty and one of $X^{IN}$ or $X^{OUT}$ is non-empty.

Now we claim that for each vertex $u$ in $ X^{OUT}$ and each vertex $v$ in $ X^{M}$, $uv \in A(D)$. Suppose this is not the case and consider $w$ an in-neighbour of $v$ in $X$: as $u$ is an out-neighbour of $w$, there must  be an arc between $u$ and $v$. If $vu \in A(D)$ then $X \cup \{v\}$ is a weak hub, contradicting the maximality of $X$. Similarly one can prove that for each vertex  $v$ in $X^{IN}$ and for each vertex $u$ in $X^{M}$, $uv \in A(D)$. 

There cannot be any arc $uv$ with $u \in X^{M} \cup X^{IN}$ and  $v \in X^{NO}$ as any such $u$ has an out-neighbour in $X$ and this would create an induce $S^+_2$. 
For a similar reason, there cannot be an arc $uv$ with $u \in X^{NO}$ and $v \in X^{M} \cup X^{OUT}$. 
In particular we prove that there are no arcs between $X^M$ and $X^{NO}$. 
But now since $X^M$ is not empty, consider a vertex $w\in X^M$ and remember that it out-dominates every vertex in $X^{IN}$ and in-dominates every vertex in $X^{OUT}$. 
Now, since $D$ is 
locally semicomplete, it does not contain $S^-_2$ nor $S^+_2$, so
there cannot be an arc from $X^{OUT}$ to $X^{NO}$ or from $X^{NO}$ to $X^{IN}$. 
Therefore there is no arc between $X^{NO}$ and  $X \cup X^{IN} \cup X^{OUT} \cup X^{M}$, and since $D$ is connected it implies that $X^{NO} = \emptyset$.

Thus $V = X \cup X^{IN} \cup X^{OUT} \cup X^{M}$ and taking $E = X^{IN}$, $F = X$, $G = X^{OUT}$ and $H = X^{M}$, we can verify that at least three of them are non-empty ($F$, $H$ and one of $E$ and $G$) and thus that they all are semicomplete (as the in-neighbourhood and the out-neighbourhood of any vertex are semicomplete) and satisfy the output of the theorem. 
\end{proof}

\subsection*{Two short applications of Theorem \ref{thm:k12k21}}
Let us mention here some applications of the previous theorem. 

A {\em 2-king} in a digraph is a vertex that can reach every other vertex by a directed path of length at most $2$. In \cite{WYW} it is proved that a locally semicomplete digraph that is not a blow-up of a round oriented graph (that is, not in the second case of Theorem \ref{thm:k12k21}) has a $2$-king. As we show now, this is a direct consequence of Theorem \ref{thm:2dic}. Indeed, when $D$ is semicomplete, it is easy to see that a vertex $x$ of maximum out-degree satisfies that every vertex in $x^-$ is either in $x^+$ or has an in-neighbour in $x^+$, and therefore $x$ is a $2$-king. And if $D$ is described by the third case, we distinguish two cases. 
\begin{itemize}
\item If $H$ is non empty, then we claim that a vertex $x$ in $G$ that is a $2$-king in $D[G]$ (which exists by the previous argument) is a $2$-king in $D$. Indeed, every vertex in $H$ is in $x^+$, and since there is all arcs from $H$ to $E$, there is a directed path of length at most $2$ from $x$ to each vertex of $E$. Finally, since $x$ has at least one out-neighbor in $E$ and there is all arcs from $E$ to $F$, there is a directed path of length $2$ from $x$ to each vertex in $F$. 
\item If $H$ is empty, then $F$ is not and in that case we take any vertex that is a $2$-king in $E$, and it is clearly a $2$-king of $D$. 
\end{itemize}

Another topic is pancyclicity, that is the property that a digraph contains a directed cycle for all possible lengths between $3$ and the number of vertices. We don't write the proof here but Theorem \ref{thm:k12k21} implies the result of \cite{BGGV} characterizing pancyclicity for locally semicomplete digaphs, since again the only non easy case is when the digraph is neither semicomplete nor the blowup of  a round digraph, in which case it follows without too much effort because of the simplicity of this third case compared to the one of Theorem 3.2 in \cite{BGGV}.

\section{Perspectives}
In the context of Conjecture  \ref{conj:oriented}, it would be interesting to prove that $\{\ovlra{K_2}, H, S^+_2\}$ is heroic for every hero $H$. 
We prove it in the case where the hero consists of a $C_3$ plus a vertex dominating it. In order to extend this partial result, one idea could be to use the structure theorem for heroes of Berger et al. mentioned in the introduction. Let us now state this theorem.

Given two oriented graphs $H_1$ and $H_2$ on disjoint sets of vertices, we denote by $H_1\Rightarrow H_2$ the digraph obtained from disjoint union of $H_1$ and $H_2$ by adding an arc from each vertex of $H_1$ to each vertex of $H_2$. 
Given three oriented graphs $H_1$, $H_2$, $H_3$, we denote by $C_3(D_1, D_2, D_3)$ the oriented graph obtained from disjoint copies of $H_1$, $H_2$ and $H_3$ after adding all arcs from $H_1$ to $H_2$, all arcs from $H_2$ to $H_3$ and all arcs from $H_3$ to $H_1$.  

\begin{theorem}[\cite{hero}]\label{thm:hero}
A tournament is a hero if and only if it can be constructed by the following inductive rules:
\begin{itemize}
	\item $K_1$ is a hero.
	\item If $H_1$ and $H_2$ are heroes, then $H_1 \Rightarrow H_2$ is also a hero.
	\item If $H$ is a hero, then for every $k \ge 1$, the tournaments $C_3(H, TT_k, K_1)$ and $C_3(TT_k, H, K_1)$ are both heroes. 
\end{itemize}

\end{theorem}

In the light of this theorem, a first step to prove that $\{\ovlra{K_2}, H, S^+_2\}$ is heroic for every hero $H$ would be to show that, given two heroes $H_1$ and $H_2$, if $\{\ovlra{K_2}, H_i, \SD\}$ is heroic for $i=1,2$, then $\{\ovlra{K_2}, H_1\Ra H_2, \SD\}$ is also heroic. 
The case $H_1=K_1$ is solved in the prepublication \cite{RS21} mentioned in the introduction, but we think that  a decomposition theorem for locally out-tournament oriented graphs, in the spirit of Theorem \ref{thm:k12k21} would  be of great help for the general case


Finally, we wonder whether Theorem \ref{thm:k12k21} can be applied to prove some open problems about locally semicomplete digraphs or locally tournament oriented graphs. One such problem is the famous Second Neighbourhood Conjecture due to Seymour, which is known for tournaments but not for locally tournament oriented graphs. The first two cases of Corollary \ref{thm:k12k21} can be dealt with without too much difficulty, but we were alas not able to deal with the last one.
\medskip 

{\bf Acknowledgment} This project was financed by the ANR projects DISTANCIA (ANR-17-CE40-0015),  HOSIGRA (ANR-17-CE40-0022) and ALGORIDAM (ANR-19-CE48-0016).
We would also like to thanks Reza Naserasr  for fruitful discussions.


\begin{thebibliography}{15}

\bibitem{digraph}
P.~Aboulker,  J. Bang-Jensen, N. Bousquet, P. Charbit, F. Havet,  F. Maffray,  and J. Zamora.
\newblock $\chi$-bounded families of oriented graphs,
\newblock {\em   Journal of Graph Theory}, 89, 3: 304--326, 2018.

\bibitem{ACN21}
P.~Aboulker, P. Charbit, R.~Naserasr.
\newblock Extension of the \Gya-Sumner conjecture to digraphs
\newblock submitted, arXiv:2009.13319.

\bibitem{ACL19}
P. Aboulker, N. Cohen, W. Lochet, F. Havet, P. Mourra, S. Thomassé
\newblock Subdivisions  in digraphs of large out-degree or large dichromatic number.
\newblock{Electronic journal of Combinatorics}, Vol. 6, 3, 2019

\bibitem{B90} J. Bang-Jensen. 
\newblock Locally semicomplete digraphs: a generalization of tournaments.
\newblock {\em Journal of Graph Theory,  14(3):371–390, 1990.}

\bibitem{BGGV}
J.Bang-Jensen, Y. Guo, G. Gutin, L. Volkmann,
\newblock A classification of locally semicomplete digraphs
\newblock {\em Discrete Mathematicss 167–168:101-114,1997.}

\bibitem{BG01}
J. Bang-Jensen, G. Gutin. 
\newblock Digraphs. Theory, algorithms and applications. Springer-Verlag London, Ltd., London, 2001.

\bibitem{BGCD}
J. Bang-Jensen, G. Gutin. 
\newblock Classes of Directed Graphs. Springer-Verlag London, Ltd., London, 2018.



\bibitem{hero}
E~Berger, K.~Choromanski, M.~Chudnovsky, J.~Fox, M.~Loebl, A.~Scott, P.~Seymour and S.~Thomass\'e,
\newblock{Tournaments and colouring},
\newblock{Journal of Combinatorial Theory, Series B}, 112:1--17, 2015.

\bibitem{BHL18}
J. Bensmail, A. Harutyunyan, N. K. Le,
\newblock{List colouring digraphs}
\newblock{Journal of Graph Thoery}, 87:492-508, 2018.

\bibitem{Perfect}
M.~Chudnovsky and N.~Robertson and P.~Seymour and R.~ Thomas,
 \newblock{The Strong Perfect Graph Theorem},
 \newblock{\em Annals of Mathematics}, 1:51-229.



\bibitem{CSS19}
M.~Chudnovsky, A.~Scott and P.~Seymour
\newblock Induced subgraphs of graphs with large chromatic number XI. Orientations. .
\newblock {\em European Journal of Combinatorics}, 76:53--61, 2019. 

\bibitem{Erd59}
P.~Erd{\H{o}}s.
\newblock Graph theory and probability.
\newblock {\em Canad. J. Math.}, 11:34--38, 1959.

\bibitem{ErHa76}
P.~Erd{\H{o}}s and A. Hajnal.
\newblock On the chromatic  number of infinite graphs.
\newblock In {\it Theory of Graphs}, {\it Proceedings of the 1966 Colloquium at Tihany}, pages 83--98,Academic Press,  1976.

\bibitem{GSS20}
L. Gishboliner, R. Steiner, T. Szabo.
\newblock Dichromatic number and forced subdivisions.
\newblock{\em Submitted}, 2020


\bibitem{Gya87}
A.~Gy\'arf\'as.
\newblock Problems from the world surrounding perfect graphs.
\newblock {\em Zastowania Matematyki Applicationes Mathematicae}, XIX:413--441,
1987.

\bibitem{HKMR}  
A. Harutyunyan, P.M. Kayll, B. Mohar, L. Rafferty, 
\newblock Uniquely D-colourable digraphs with large girth, 
\newblock {\em The Canadian Journal of Mathematics}, 64:1310-1328, 2012. 


\bibitem{girth}
A. Harutyunyan and B. Mohar, Two  results  on  the  digraph  chromatic  number, Discrete Math., 312:1823–1826, 2012.

\bibitem{HLNT}
A. Harutyunyan, T-N. Le, A. Newman, S Thomassé, 
\newblock colouring dense digraphs.
\newblock  {\em Combinatorica}, 39:1021--1053, 2019.
 \bibitem{H17}
W. Hochst\"attler
\newblock A flow theory for the dichromatic number
\newblock {\em European Journal of Combinatorics}, 66:160--167, 2017. 


\bibitem{HK15}
R. Hoshino and K. Kawarabayashi 
\newblock The edge density of critical digraphs.
\newblock  {\em Combinatorica}, 35:619--631, 2015.
 
\bibitem{H92} J. Huang. Tournament-like Oriented Graphs. PhD thesis, CAN, 1992.


 \bibitem{KiRo96}
 H.~A. Kierstead and V.~R{\"o}dl.
 \newblock Applications of hypergraph colouring to colouring graphs that do not
 induce certain trees.
 \newblock {\em Discrete Math.}, 150:187--193, 1996.

\bibitem{KiTr92}
H.~A. Kierstead and W. T.~Trotter.
\newblock colourfull induced subgraphs,
\newblock {\em Discrete Math.}, 101:165--169, 1992.

\bibitem{KS20}
A.V. Kostochka and M. Stiebitz 
\newblock The Minimum Number of Edges in 4-Critical Digraphs of Given Order. 
\newblock {\em Graphs and Combinatorics}, 36:703--718, 2020. 

\bibitem{LM17} Z. Li and B. Mohar. Planar digraphs of digirth four are 2-colourable.
\newblock{\it SIAM J. Discrete Math.}, 31:2201--2205, 2017.

\bibitem{LZM08} R. Li, X. Zhang, and W. Meng.
\newblock  A sufficient condition for a digraph to be positive round.
\newblock {\em Optimization, 57:345–352, 2008.}

\bibitem{L11} G. Liu, Digraphs constructed by iterated substitution from a base set, Junior paper, Princeton, 2011.

\bibitem{NL82}
V. Neumann-Lara.
\newblock The dichromatic number of a digraph. 
\newblock {\em Journal of Combinatorial Theory, Series B}, 33:265--270, 1982.

\bibitem{R13}
A.A. Razborov,
\newblock On the Caccetta–Häggkvist Conjecture with Forbidden Subgraphs. 
\newblock{\em J. Graph Theory, 74: 236-248}, 2013.


\bibitem{SS20}
A. Scott, P. Seymour, 
\newblock A survey of $\chi$-boundedness, 
\newblock {\em Journal of Graph Theory} 95:473-504, 2020.


\bibitem{SS16-1}
A. Scott and P. Seymour
\newblock Induced subgraphs of graphs with large chromatic number. I. Odd holes
\newblock {\em Journal of Combinatorial Theory, Series B},121:68-84, 2016.
	
  
\bibitem{Stearns}
\newblock R. Stearns, The voting problem, 
\newblock {\em Amer. Math. Monthly 66 (1959), 761--763.}

\bibitem{RS21}
R. Steiner
\newblock On coloring digraphs with forbidden induced subgraphs
\newblock {\em arXiv Preprint 2103.04191, 2021    }

\bibitem{S20}
R. Steiner. 
\newblock A note on coloring digraphs of large girth, 
\newblock{\em Discrete Applied Mathematics}, 287:62-64, 2020.

\bibitem{Sum81}
D.~P. Sumner.
\newblock Subtrees of a graph and the chromatic number.
\newblock {\em The theory and applications of graphs (Kalamazoo, Mich.,
  1980)}, pages 557--576. Wiley, New York, 1981.
  
\bibitem{WYW} R. Wang, A. Yang, and S. Wang. 
\newblock Kings in locally semicomplete digraphs.
\newblock {\em Journal of Graph Theory, 63(4):279–287, 2010.}

 
\end{thebibliography}
\end{document}